\documentclass[11pt]{article}
\catcode`\@=11 \catcode`\@=12

\usepackage{amsmath}
\usepackage{amssymb}
\usepackage{amsthm}
\usepackage{oldgerm}

\usepackage[pdftex,colorlinks,urlcolor=blue,pdfstartview=FitH]{hyperref}

\newcommand{\Rm}{\mathbb{R}}

\newcommand{\mL}{\mathcal{L}}

\newcommand{\mC}{\ensuremath{\mathcal{C}}}

\newcommand{\mS}{\ensuremath{\mathcal{S}}}
\newcommand{\mF}{\ensuremath{\mathcal{F}}}

\newcommand{\Qm}{\ensuremath{\mathbb{Q}}}
\newcommand{\Nm}{\ensuremath{\mathbb{N}}}

\newcommand{\mM}{\ensuremath{\mathcal{M}}}

\newcommand{\mY}{\ensuremath{\mathcal{Y}}}
\newcommand{\mA}{\ensuremath{\mathcal{A}}}

\newcommand{\mB}{\ensuremath{\mathcal{B}}}

\newtheorem{lem}{Lemma}
\newtheorem{thm}{Theorem}

\newtheorem{cor}[lem]{Corollary}
\newtheorem{prop}[lem]{Proposition}
\newtheorem{defn}[lem]{Definition}

\def\proof {\noindent{\sc{Proof. }}}
\def\qed {\mbox{}\hfill {\small \fbox{}} \\}
\def\lto{\longrightarrow}
\def\lmto{\longmapsto}

\def\leq{\leqslant}
\def\geq{\geqslant}

\title{Some remarks on the continuity equation}

\author{Patrick Bernard\footnote{membre de l'IUF} }

\date{}

\begin{document}

\maketitle

\begin{small}
\begin{center}
------\\
\noindent
Patrick Bernard,
 Universit\'e  Paris-Dauphine,\\
CEREMADE, UMR CNRS 7534\\
Pl. du Mar\'echal de Lattre de Tassigny\\
75775 Paris Cedex 16,
France\\
\texttt{patrick.bernard@ceremade.dauphine.fr}\\
------
\end{center}

\vspace{1cm}

\end{small}

This text is the act of a talk given
november 18 2008 at the seminar  PDE
of Ecole Polytechnique. The text is not completely faithfull
to the oral exposition for I have taken this opportunity 
to present the proofs of some results that are 
not easy to find in the literature. On the other hand, 
I have been less precise  on the material 
for which I found good references. Most of the novelties
presented here come from a joined work with Luigi Ambrosio.

\section{Introduction}
We consider a Borel vector-field
$$
V(t,x):]0,T[\times \Rm^d\lto \Rm^d,
$$
and the associated equations 
\begin{equation}\tag{ODE}\label{ODE}
\dot \gamma(t)=V(t,\gamma(t))
\end{equation}
and  (with the notations  $V_t(x)=V(t,x)$)
\begin{equation}\tag{PDE}\label{PDE}
\partial_t \mu_t + \text div(V_t \mu_t)=0.
\end{equation}
It is important to notice that $V(t,x)$ is a well-defined
function, and not an equivalence class of functions.
In order to avoid some technicalities we assume the bound
\begin{equation}\tag{B}\label{B}
\|V\|_c:=\int_0^T \|V_t\|_{\infty}dt <\infty.
\end{equation}
Here $\|V_t\|_{\infty}$ is defined as the supremum of $\|V(t,x)\|$.
A solution of (\ref{ODE}) is an absolutely continuous curve $\gamma(t)$
such that $\dot \gamma(t)=V(t,\gamma(t))$ almost everywhere on
$[0,T]$.

We  consider solutions of (\ref{PDE}) in the class
$\mM(\Rm^d)$
of bounded signed measures. It is necessary here to settle
a couple of notations.
We define the Banach space $C_0(\Rm^d)$ as the 
set of continuous functions which converge
to zero at infinity. It is endowed with the uniform 
norm. The space $\mM(\Rm^d)$ is the dual of $C_0(\Rm^d)$,
we endow it with the weak-$*$ topology, that we will 
simply call the weak topology.
We denote by $\mM^+(\Rm^d)$ and $\mM^{1+}(\Rm^d)$
the spaces of non-negative and probability Borel measures.
 Given a signed  measure
 $\mu\in \mM(\Rm^d)$, we denote by $|\mu|\in \mM^+(\Rm^d)$
 its total variation.
 The quantity $\|\mu\|:= |\mu|(\Rm^d)$ defines a norm on 
 $\mM(\Rm^d)$, which coincides with the dual norm.

A solution of (\ref{PDE}) is a weakly continuous
curve $\mu_t\in C([0,T],\mM(\Rm^d))$ such that,
for each compactly supported smooth function $u$
on $\Rm^d$, the function 
$t\lmto \int_{\Rm^d}u(x)d\mu_t(x)$ is absolutely continuous
with derivative given by
\begin{equation}\tag{D}\label{D}
\left(\int_{\Rm^d} u(x)d\mu_t(x)\right)'=\int _{\Rm^d}du_x(V(t,x))d\mu_t(x).
\end{equation}
This relation then holds for each $C^1$ function $u$ which
is bounded and Lipschitz.
That this definition is equivalent to the genuine definition 
in the sense of distributions
(and in particular, that weak continuity is in fact a consequence
of being a solution) is explained, for example, in  
\cite{AGS}, Chapter 8.
Note that, for non-negative solutions, the norm 
$\|\mu_t\|=\int 1d\mu_t$ is preserved. For signed solutions,
however, this norm is not necessarily continuous and may not 
be bounded. We will restrict our attention 
to norm-bounded solutions (those for which the function 
$t\lmto \|\mu_t\|$
is bounded).

In the good cases, (\ref{ODE}) can be solved by a flow:

\begin{defn}\label{flow}
The Borel map
$$X(t,s,x):[0,T]\times [0,T]\times \Rm^d\lto \Rm^d
$$ 
is called a flow of solutions of (\ref{ODE})
if, for each fixed $t$ and $x$, the curve $s\lmto X(t,s,x)$
is the only solution $\gamma(s)$ of (\ref{ODE})
which satisfies the initial condition $\gamma(t)=x$.
 The maps $X_t^s:\Rm^d\lto \Rm^d$
defined by $X_t^s(x):=X(t,s,x)$
then satisfy the Markov property
$$
X_{t_1}^{t_2}\circ X_{t_0}^{t_1}=X_{t_0}^{t_2}.
$$
\end{defn}

It follows from Proposition \ref{mf} below that a
flow of solutions of (\ref{ODE})
exists if and only if, for each $S\in [0,T]$ and 
$x\in \Rm^d$, the Cauchy problem 
consisting of solving (\ref{ODE}) with the initial data
$\gamma(S)=x$ has one and only one solution.
If $X(t,s,x)$ is the flow of solutions of (\ref{ODE}), 
then it is easy to see that, for each given $S\in [0,T]$
and $\mu\in \mM$, the expression
$$
\mu_t:= (X_S^t)_{\sharp} \mu
$$
defines a solution of (\ref{PDE}).
We say that the flow $X$ uniquely solves (\ref{PDE})
if this is the only norm-bounded solution of (\ref{PDE})
fulfilling the given initial value.

When the vector-field $V(t,x)$ is smooth,  there exists
a smooth flow $X(t,s,x)$
which uniquely solves (\ref{ODE}) and (\ref{PDE}).
One of our goals in the present text is to present
more general classes 
 of vector-fields 
for which both (\ref{PDE}) and (\ref{ODE}) are uniquely solved
by a flow.
 Before this, we settle some measurability issues in 
Section \ref{mes}, and then describe, 
following Ambrosio, Gigli and Savar\'e, some 
important relations between
(\ref{ODE}) and (\ref{PDE}) which hold in full generality.
They imply in particular that a flow solving (\ref{ODE})
always solves (\ref{PDE}) uniquely in the class of non-negative
measures. The situation is more intricate in the class of signed measures.
In order to understand this fact, it is important to realize
that the positive part of a signed solution 
is not necessarily a solution in general.
For example, assume that there exists two solutions 
$x(t)$ and $y(t)$ of (\ref{ODE}) and a time $S\in ]0,T[$
such that $x(S)=y(S)$. Then we can define the signed solution 
$\mu_t$ by $\mu_t=0$ for $t\leq S$ and 
$\mu_t=\delta_{x(t)}-\delta_{y(t)}$ for $t\geq S$.
It is not hard to see that this is a solution of (\ref{PDE}),
but that the positive part is not.
This example also illustrates the non-continuity of the norm
$t\lmto \|\mu_t\|$ for signed solutions.

In order to study the existence of flows, we then focus our attention
to vector-fields which are continuous in the space variable.
We denote by $\mC([0,T]\times \Rm^d)$ (or simply $\mC$) the set of 
Borel vector-fields
$V(t,x)$  such that, for each $t$,
the map $V_t:x\lmto V(t,x)$ is continuous, and such that, 
in addition, the estimate (\ref{B}) holds. The quantity 
$\|V\|_c$ defined by $\|V\|_c:= \int_0^T \|V_t\|_{\infty}dt$
is a norm on $\mC$, and $\mC$ endowed with this norm
is a Banach space.
If $V\in \mC$, it is well-known that, for each $S$ and $x$, there
exists a solution $\gamma$  to (\ref{ODE}) satisfying $\gamma(S)=x$.
The existence of a flow of solutions of (\ref{ODE})
is then  equivalent to the uniqueness for each Cauchy data.
In Section \ref{generic}, we will prove:

\begin{thm}\label{gen}
The set of vector-fields $V\in \mC$ for which both (\ref{ODE})
and (\ref{PDE}) are uniquely  solved by a flow is generic in $\mC$
in the sense of Baire.
\end{thm}

This result is rather easy, but
we do not know any reference, and will therefore provide a 
complete proof.
A different but similar genericity
result is proved in \cite{L}.
Our proof is very different from the one in \cite{L},
but a similar one could possibly be used.

Next we try to derive the existence of a flow from 
regularity estimates.
We recall that a modulus of continuity is a continuous
non-decreasing function $\rho:[0,1)\lto [0,\infty)$, such that
$\rho(0)=0$. A modulus of continuity $\rho$ is said to be Osgood if
$$
\int_0^1 \frac{1}{\rho(s)} ds=+\infty.
$$
We will always extend the moduli of continuity to $[1,\infty)$ by
$\rho=\infty$. Typical examples of Osgood moduli of continuity are
$\rho(s)=s$ and $\rho(s)= s(1-\ln(s))$. Note that the moduli
$\rho(s)=s^\alpha$, $\alpha\in (0,1)$, are not Osgood.

It is known that (\ref{ODE}) is solved by a unique flow
(which is a flow of homeomorphisms) provided
 there exists an Osgood
modulus of continuity $\rho$ and $C(t)\in L^1(0,T)$ such that
\begin{equation}\tag{O}\label{O}
| V(t,x)-V(t,y)|\leq
C(t)\rho(|x-y|)
\end{equation}
for all $x,\,y\in\Rm^d$, and all $t\in ]0,T[$.
We do not know if, in general, the existence of a flow of 
homeomorphisms  solving (\ref{ODE}) implies that 
this flow also uniquely solves (\ref{PDE})
(although we know that it uniquely solves (\ref{PDE})
in the class of non-negative measures, see Section \ref{superposition}).
If $V\in \mC$ satisfies (\ref{O}), then this is true:

\begin{thm}\label{os}
If $V\in \mC$ satisfies (\ref{O}), then there exists a flow of
homeomorphisms uniquely solving (\ref{ODE}) and (\ref{PDE}).
\end{thm}

This result was proved in \cite{AB}. 
It had been proved earlier in \cite{BC} in the case where 
$\rho(s)=s(1-\ln(s))$ and where $V$ is incompressible.
The method was strikingly different.
The result  can be considered standard
in the case $\rho(s)=s$.
We  give some indications of proof in Section \ref{osgood}.

\section{Some measurability issues}\label{mes}
\begin{prop}\label{Borel}
If $V$ is a Borel vector-field, then the set of solutions 
of (\ref{ODE}) is Borel in $C([0,T],\Rm^d)$.
\end{prop}
\proof
The curve $\gamma\in C([0,T],\Rm^d)$ is a solution if and only if
$$
\gamma(t)-\gamma(0)=\int_0^t V(s,\gamma(s))ds
$$
for each $t\in \Qm\cap [0,T]$.
So in order to prove the Proposition, it is enough
to prove that, for a given $t\in [0,T]$, the map 
$\gamma\lmto \int_0^t V(s,\gamma(s))ds$
is Borel.
We claim that, for each non-negative Borel function 
$f(s,x):[0,t]\times \Rm^d\lto [0,\infty)$, the map
$$\gamma\lmto \int_0^tf(s,\gamma(s))ds$$
is Borel. The claim clearly implies the desired result.

We prove the claim by a monotone class argument.
Let $\mF$ be the set of functions which satisfy 
the desired conclusion.
We observe that $\mF$ is stable under addition, 
multiplication by a non-negative real number, and monotone 
convergence. Moreover, $\mF$ obviously contains bounded continuous
functions.
By standard monotone class arguments, we conclude that $\mF$
contains all non-negative Borel functions
(see \textit{e. g.} \cite{B}, Lemma 39). 
\qed

\begin{prop}\label{mf}
Let $B$ be a Borel set in $[0,T]\times \Rm^d$
such that, for each $(S,x)\in B$,
there 
exists one and only one solution $\gamma(t)$ of (\ref{ODE})
which satisfies $\gamma(S)=x$.
Then the flow map $X:[0,T]\times [0,T] \times \Rm^d\lto \Rm^d$
is well-defined and Borel on the set 
$$
\tilde B:=\{(S,t,x)\in [0,T]\times [0,T]\times \Rm^d
\quad |\quad (S,x)\in B\}.
$$
\end{prop}
\proof
Let $\mA\subset C([0,T],\Rm^d)$ be the (Borel) set of solutions
of (\ref{ODE}).
The map 
$$ev:[0,T]\times \mA\lto [0,T]\times \Rm^d
$$
given by $ev(t,\gamma)=(t,\gamma(t))$ is continuous hence Borel.
As a consequence, the set $\mA_R:=ev^{-1}(R)$ is Borel.
Our hypothesis is precisely that the map $ev$
is a bijection between $\mA_R$ and $R$.
By (non trivial) general results  
of measure theory (see \cite{K} or \cite{P}, Theorem 3.9), the inverse map 
$ev^{-1}$ is then Borel.
Now the conclusion follows from the  formula
$$
X(S,t,x)=\pi\circ ev(t,\pi\circ ev^{-1}(S,x)),
$$
where we have denoted by the same letter $\pi$ two
different projections consisting in forgetting the 
time.
\qed

\section{The uniqueness question for probability measures}
\label{superposition}
We describe several links between 
(\ref{ODE}) and (\ref{PDE}) which hold when 
$V$ is only Borel. The content of this section is due to 
Ambrosio, Gigli
and Savar\'e (see \cite{AGS}, chapter 8),
but is also closely related to earlier works of 
Smirnov (\cite{S}) as will be made more apparent in section \ref{osgood}.

Let us first recall that every solution of (\ref{ODE})
can be seen as a non-negative solution of (\ref{PDE}).
Indeed, if $\gamma(t)$ solves (\ref{ODE}), then
the measures $\mu_t:= \delta_{\gamma(t)}$ solves (\ref{PDE}).
We call the solutions of (PDE) which can be obtained this way 
elementary.
The following statement, established by Ambrosio, Gigli
and Savar\'e in the line of anterior works of Smirnov,
roughly states that the set of solutions of (\ref{PDE})
in the class of probability measures coincides with the closed
convex envelop of the set of elementary solutions.

\begin{thm}\label{AGS}
Let $\mu_t\in C([0,T],\mM^{1+})$ be a solution of (\ref{PDE}).
Then there exists a Borel probability measure 
$\nu$ on $C([0,T],\Rm^d)$ with the following properties:
\begin{enumerate}
\item $\nu$ is concentrated on the Borel set of solutions of 
(\ref{ODE}), 
\item $(ev_t)_{\sharp}\nu=\mu_t$ for each $t\in [0,T]$, where 
$ev_t:C([0,T],\Rm^d)\lto \Rm^d$ is the evaluation map
$\gamma\lmto \gamma(t)$.
\end{enumerate}
\end{thm}

This theorem can be equivalently stated as follows:

\begin{thm}\label{process}
Let $\mu_t\in C([0,T],\mM^{1+})$ be a solution of (\ref{PDE}).
Then there exists a stochastic process $Z(t,\omega)$
such that :
\begin{enumerate}
\item almost each sample path $t\lmto Z(t,\omega)$ is a solution of (\ref{ODE})
\item The law of the random variable $\omega\lmto Z(t,\omega)$
is $\mu_t$.
\end{enumerate}
\end{thm}

These results, called the superposition principle,
 sum up what is known in general 
concerning the relations between (\ref{ODE}) and(\ref{PDE}).
Excellent proofs can be found in 
\cite{AGS,Am:cetraro,AC}, see also \cite{B} for another approach,
and \cite{maniglia,BB} for related material.
The following corollary is obvious:

\begin{cor}\label{exist}
Let $\mu_t\in C([0,T],\mM^{1+})$ be a solution of (\ref{PDE}).
Then, for each $S\in [0,T]$ and for $\mu_S$-almost every
point $x\in \Rm^d$, there exists a solution $\gamma(t)$
of (\ref{ODE}) such that $\gamma(S)=x$.
\end{cor}

\begin{cor}
If $V_t$ is volume-preserving for each $t$,
then, for each $S\in [0,T]$ and for almost every
point $x\in \Rm^d$, there exists a solution $\gamma(t)$
of (\ref{ODE}) such that $\gamma(S)=x$.
\end{cor}

\proof
This corollary does not immediately follow from Corollary \ref{exist}
because the Lebesgue measure is not bounded.
We denote by $\lambda$ the Lebesgue measure, and consider 
a positive and bounded function
 $v:\Rm^d\lto ]0,1]$ such that $\int vd\lambda =1$, so that 
 $\mu=f\lambda$ is a probability measure. The corollary holds
if there exists a solution $\mu_t\in C([0,T],\mM^+)$
such that $\mu_S=\mu$.
This is the content of the following Lemma:
\qed

\begin{lem}
Let $V(t,x)$ be a Borel vector-field
such that $\|V\|_c<\infty$ and $div(V_t)=0$ for each $t$.
Let $S\in [0,T]$ be a fixed  time and let $v:\Rm^d\lto [0,1]$
be an integrable function (normalized to $\int vd\lambda =1$).
Then there exists a non-negative  solution 
$\mu_t\in C([0,T],\mB)$  of (\ref{PDE})
with $\mu_S=v\lambda$ and such that,
for each time $t\in [0,T]$, we have $\mu_t=v_t\lambda$
for some integrable function $v_t:\Rm^d\lto [0,1]$.
\end{lem}
\proof
We mollify $V$  by
\begin{equation}\tag{M}\label{M}
W^n(t,x):= n^d\int V(t,y)g(n(x-y))dy,
\end{equation}
where $g$ is a compactly supported smooth kernel.
We have $W^n\lto W$ in $L^1_{loc}$, and $W^n$ satisfy
the estimate (\ref{O}) with a Lipschitz modulus.
As a consequence, for each $n$, there exists a 
flow of homeomorphisms solving (\ref{ODE}), and therefore
there exists a  Borel  function 
$v^n(t,x):[0,T]\times \Rm^d\lto [0,1]$ such that 
$\mu^n_t:=v^n_t\lambda$ is a solution of (\ref{PDE})
(with the vector-field $W^n$)
and such that $\mu^n_S=v\lambda$.
Observe also that $\|W^n\|_c\leq \|V^n\|_c$, from which follows,
using (\ref{E}) in the appendix, that the sequence $\mu^n_t$
is equicontinuous.
As a consequence, we can assume  that $\mu^n_t$
converges uniformly  to  a limit $\mu_t\in \mC([0,T],\mB)$.
Note that $\mu_S=v\lambda$.
For each fixed $t$, 
the measure $\mu^n_t$ has a density $v^n_t$
with values in $[0,1]$
hence  the limit $\mu_t$ has a density $v_t$ with 
values in $[0,1]$, and $v_t^n\lto v_t$ weakly-$*$
in $L^{\infty}(\Rm^d)$.
We have to prove that $\mu_t$ solves (\ref{PDE}).

We have
\begin{equation}\label{eqlem}
\int_0^Tf'(t)\int u(x)d\mu^n_t(x)   +
f(t) \int du_x(W^n(t,x))d\mu^n_t(x)\quad dt=0
\end{equation}
for each compactly supported smooth function $u$ on $\Rm^d$
and each smooth compactly supported  function $f$ on $]0,T[$.
For each fixed $t$ the functions
$du_x(W^n_t(x))$ strongly converge 
to $du_x(V_t(x))$ in $L^1$.
Since in addition $v^n_t$ converges to $v_t$
weakly-$*$ in $L^{\infty}(\Rm^d)$, we have:
$$
\int du_x(W^n(t,x))v^n_t(x)d\lambda(x)
\lto 
\int du_x(V(t,x))v_t(x)d\lambda(x).
$$
By the dominated convergence theorem
(using that $\|W^n\|_c$ is bounded), we can pass to the limit
in (\ref{eqlem}) and get 
$$
\int_0^Tf'(t)\int u(x)d\mu_t(x)   +
f(t) \int du_x(W(t,x))d\mu_t(x)\quad dt=0
$$
which says exactly 
that the measures $\mu_t$
solve (\ref{PDE}).
\qed

On the side of uniqueness, we have:

\begin{cor}
Let $S\in [0,T]$ be fixed, and let $B\in \Rm^d$ be a Borel set
such that, for each $x\in B$, there exists at most one solution
$\gamma(t)$ of (\ref{ODE}) satisfying $\gamma(S)=B$.
Then, if $\mu$ is a probability measure concentrated 
on $B$, there exists at most one solution 
$\mu_t\in C([0,T],\mM^+)$ of (\ref{PDE})
satisfying $\mu_S=\mu$.
\end{cor}
Note that, in general, there may exist other solutions
in $ C([0,T],\mM)$.

\proof
Let $\mu_t$ and $\tilde \mu_t$ be two solutions in 
$C([0,T],\mM^+)$ satisfying $\tilde \mu_S=\mu=\mu_S$.
Let $\nu$ and $\tilde \nu$ be the decompositions given by 
Theorem  \ref{AGS}.
We claim that $\nu=\tilde \nu$, and therefore
that $\mu_t=\tilde \mu_t$ for each $t$.
In order to prove the claim, 
we consider the Borel subset $Q\subset C([0,T],\Rm^d)$
formed by  solutions $\gamma$ of (\ref{ODE}) 
which satisfy $\gamma(S)\in B$.
That this set is Borel follows from Proposition \ref{Borel}.
Our hypothesis is that the restriction to $Q$ of the 
evaluation map  $ev_S$ is one-to-one.
As a consequence (by Theorem 3.9 in \cite{P})
the image $B'=ev_S(Q)$ is Borel and the inverse map 
$ev_S^{-1}:B'\lto Q$ is Borel.
By Theorem \ref{AGS}, we have $\nu(Q)=1$,
so the measure $\nu$ can be identified with 
its restriction to $Q$. 
Since $(ev_S)_{\sharp}\nu=\mu$, we have 
$\mu(B')=\nu(Q)=1$, so that  the measure $\mu$
coincides with its restriction to $B'$.
As a consequence, we obtain
$\nu=(ev_S^{-1})_{\sharp}\mu$.
Similarly, we have $\nu'=(ev_S^{-1})_{\sharp}\mu$.
\qed

\begin{cor}
If (\ref{ODE}) is uniquely solved by the flow $X$,
then, for each $S\in [0,T]$ and each probability measure
$\mu$, the curve
$\mu_t=(X_S^t)_{\sharp}\mu$ is the only solution 
of (\ref{PDE}) in $C([0,T],\mM^+)$ which satisfies 
$\mu_S=\mu$.
\end{cor}

There may exist other norm-bounded solutions in 
$C([0,T],\mM)$.
It is not easy to give good extensions of the theory
presented in this section for the case of signed measures.
I shall present some recent works in that direction in Section \ref{osgood}.

\section{Generic existence of a flow}\label{generic}
In this section, we prove Theorem \ref{gen}.
In order to prove that both (\ref{ODE})
and (\ref{PDE}) are solved by a flow for a given vector-field $V$,
it is enough to prove that any norm-bounded solution $\mu_t$
of (\ref{PDE}) which satisfies $\mu_0=0$ or $\mu_T=0$
must vanish identically.
Indeed, assume that there exists two   solutions $\mu_t$
and $\mu'_t$
which are not equal, and some $S\in [0,T]$
such that $\mu_S=\mu'_S$. 
Let us assume for instance that there exists
$t\geq S$ such that $\mu_t\neq \mu'_t$. Then we can define a new 
solution 
$\tilde \mu_t$ by $\tilde \mu_t=0$ for $t\leq S$
and $\tilde \mu_t=\mu_t-\mu'_t$ for $t\geq S$. The solution $\tilde \mu_t$ vanishes
at $t=0$ and it is not identically zero.

So what we have to prove now is that, for a generic vector-field
$V$, there is no  non-trivial solution of (\ref{PDE})
satisfying $\mu_0=0$ (the analogous statement for
$\mu_T$ is similar).

Let us define the set-valued mapping $\mS$ which, 
to each  vector-field $V\in  \mC$, associates the subset 
$\mS(V)\subset C([0,T],\mB)$
formed by those solutions $\mu_t$ of (\ref{PDE}) which vanish 
at time $t=0$.
As explained in Appendix \ref{top}, we embed $C([0,T],\mB)$
into a compact metric space $\mY$, and consider
$\mS$ as a set-valued map between $\mC$ and $\mY$.
We refer to Appendix \ref{set} for the terminology on
set-valued maps.

\begin{lem}
The set-valued map $\mS$ has closed graph (or equivalently it is upper semi-continuous)
\end{lem}

\proof
Let $V^n\in\mC$ be a sequence of vector-fields converging 
to $V$ in $\mC$, and let $\mu^n_t\in C([0,T],\mB)$ 
be a sequence of solutions of (\ref{PDE}) with vector-fields
$V^n$ converging to a limit $\mu_t\in \mY$ (see Appendix \ref{top}
for the definition of $\mY$).
We have
$$
\big|\|V_t^n\|_{\infty}-\|V_t\|_{\infty}\big|\leq
\|V_t^n-V_t\|_{\infty}
$$
and therefore the functions 
$t\lmto \|V_t^n\|_{\infty}$ converge to $\|V_t\|_{\infty}$ in $L^1([0,T])$.
As a consequence, the functions $\|V^n_t\|_{\infty}$
are equi-integrable, and we conclude from
inequality  (\ref{E}) in Appendix \ref{top} that
the curves $\mu^n_t$ are equi-continuous. 
As a consequence, the limit $\mu_t$
belongs to $C([0,T],\mB)$.
We have to prove that $\mu_t\in \mS(V)$.
Recall that  $\mu^n_t\in \mS(V^n)$ if and only if the equality 
\begin{equation}\label{test}
\int_0^T\int f'(t)u(x) d\mu^n_t(x) dt+
\int_0^T\int f(t)du_x(V^n(t,x))
d\mu^n_t(x) dt=0
\end{equation}
 holds for each compactly supported  smooth function $u$ on $\Rm^d$
and each smooth function 
$f:\Rm \lto \Rm$ satisfying $f(T)=0$.
Since $\mu^n_t$ converge to $\mu_t$
in $\mY$, we have 
$$
\int_0^T\int f'(t)u(x) d\mu^n_t(x) dt
\lto 
\int_0^T\int f'(t)u(x) d\mu^n_t(x) dt
$$
and 
$$\
\int_0^T\int f(t)du_x(V(t,x))
d\mu^n_t(x) dt
\lto
\int_0^T\int f(t)du_x(V(t,x))
d\mu_t(x) dt.
$$
On the other hand,
we have 
$$
\left|
\int_0^T\int f(t)du_x(V(t,x))
d\mu^n_t(x) dt
-
\int_0^T\int f(t)du_x(V^n(t,x))
d\mu^n_t(x) dt
\right|
\leq
C\|V_n-V\|_c\lto 0
$$
hence, at the limit in (\ref{test}),
we obtain
$$\int_0^T\int f'(t)u(x) d\mu_t(x) dt+
\int_0^T\int f(t)du_x(V(t,x))
d\mu_t(x) dt=0
$$
for each compactly supported smooth function $u$ on $\Rm^d$
and each smooth function 
$f:\Rm\lto \Rm$ satisfying $f(T)=0$.
This implies that $\mu_t\in \mS(V)$.
\qed
From Kuratowski Theorem ( Theorem \ref{K} in Appendix \ref{set} ),
we conclude that the set of points of continuity of $\mS$
is generic.
In order to end the proof, it is enough to see that
$\mS(V)=\{0\}$ when $V$ is a point of continuity of $\mS$.
Let $W^n(t,x)$ be the sequence of mollified approximations of 
$V$ defined in (\ref{M}).
We have $W^n\lto V$ in $\mC$. On the other hand, we have 
$\mS(W^n)=\{0\}$.
Since $V$ is a point of continuity of $\mS$, for each solution
$\mu_t\in \mS(V)$, there exists a sequence $\mu^n_t\in \mS(W^n)$
such that $\mu^n_t\lto \mu_t$ in $\mY$.
This implies 
that $\mS(V)=\{0\}.$
\qed

\section{The uniqueness question for signed measures}\label{osgood}
Let us first recall the following well-known result:

\begin{thm}\label{un}
Consider a vector-field  $V\in \mC$ which satisfies (\ref{O}).
Then there exists a unique flow  of homeomorphisms $X_s^t$
solving (\ref{ODE}).
\end{thm}

Our main issue here is to prove Theorem \ref{os},
that is to prove that the vector-field
which solves (\ref{ODE}) also uniquely solves (\ref{PDE}).
As far as non-negative solutions of (\ref{PDE}) are concerned,
this is implied by Section \ref{superposition}.
So it is natural to start with a new superposition principle
adapted to signed solutions.

Our method to de so is to consider the extended 
vector-field $\tilde V(t,x):=(1,V(t,x))$ on $[0,T]\times \Rm^d$.
Now if the measure $\mu_t$ solve (\ref{PDE}),
then setting $\tilde \mu:= dt\otimes \mu_t$ (extended by zero
outside of $]0,T[\times \Rm^d$), we have
$$
div(\tilde V\tilde \mu)=\delta_T\otimes \mu_T-\delta_0\otimes \mu_0
$$
on $\Rm^{d+1}$. 
Here the important object is the product $\tilde V\tilde \mu$
which is a vector-valued measure and even  a normal one-current
(a vector-valued measure whose divergence is also a vector-valued
measure).
From general results of Smirnov (see \cite{S}) on the decomposition
of normal one-currents,  we infer (see \cite{AB}) 
the following superposition principle for the signed solutions
of (\ref{PDE}):

\begin{thm}\label{sup}
Let $V(t,x)$ be a Borel vector-field satisfying (\ref{B}),
let $\mu_t\in C([0,T],\mM)$ be a solution of (\ref{PDE}), and
let $C(t)$ be a given positive integrable function on $[0,T]$.
Then 
 there exists 
a Borel probability measure $\nu$
on $C([0,1],\Rm^{d+1})$ such that :
\begin{enumerate}
\item \label{boundary}
$\delta_{T}\otimes \mu_{T}-\delta_0\otimes \mu_0
=(ev_1)_{\sharp}\nu-(ev_0)_{\sharp}\nu
$
\item \label{curve} $\nu$-almost each curve $\gamma(s)=(t(s),x(s))$ is one to one
and Lipschitz, it satisfies the estimate
$\int _0^1 C(t(s))|\dot t(s)| ds<\infty$,
it takes values in $[0,T]\times \Rm^d$, and solves the equation
\begin{equation}\tag{R}\label{R}
\dot x(s)=\dot t(s) V(t(s),x(s))
\end{equation}
for almost every $s$.
\end{enumerate}
\end{thm}

The equality in \ref{boundary} is global,  it is not true
in general 
that $\delta_{T}\otimes \mu_{T}=(ev_1)_{\sharp}\nu$
or $\delta_0\otimes \mu_0=(ev_0)_{\sharp}\nu$.

This superposition principe is far less appealing than
Theorem \ref{AGS}. This kind of complication seems unavoidable.
The price paid from the existence of a non-constant sign 
of the solution is the fact that the time component $t(s)$
of the curves appearing in the decomposition is not
necessarily monotone.
In order to understand the role of the equation (\ref{R}),
it is worth noticing that, if $x(t)$ is a solution 
of (\ref{ODE}) and if $t(s):[0,1]\lto [0,T]$ is any Lipschitz 
function, then the curve $\gamma(s)=(t(s),x(t(s)))$
satisfies (\ref{R}).

Conversely, if we could prove that (\ref{R}), seen as 
the  ODE
\begin{equation}\tag{R'}\label{R'}
\dot x(s)=P(s,x(s))
\end{equation}
with $P(s,x):[0,1]\times \Rm^d\lto \Rm^d$ given by
$$P(s,x):=\dot t(s) V(t(s),x)$$
   satisfies uniqueness, then we 
would conclude that
$$x(s)=X(t(s_0),t(s),x(s_0))$$
for each $s_0$ and $s$ in $[0,1]$, where $X$ is the flow solving 
(\ref{ODE}) with the vector-field $V(t,x)$.
In general, given $t(s)$, we do not see any reason why uniqueness
should hold for (\ref{R'}) even if it holds for (\ref{ODE}).
However, we have:

\begin{lem}\label{su}
Under the hypotheses of Theorem \ref{sup}, if in addition 
the vector-field $V$ satisfies (\ref{O}) with the same function 
$C(t)$ as in Theorem \ref{sup},
then $\nu$-almost
each curve $\gamma(s)=(t(s),x(s))$
satisfies 
$$x(s)=X(t(0),t(s),x(0)).$$
\end{lem}

\proof
We have 
$
|P(s,y)-P(s,x)|\leq D(s)|x-y|,
$
with $D(s)=|\dot t(s)|C(t(s))$.
For $\nu$-almost every curve $\gamma(s)=(t(s),x(s))$,
 we know that the function
$D(s)=|\dot t(s)|C(t(s))$ is integrable
(this is one of the conclusions of Theorem \ref{sup}).
In other words, the vector-field $P$ satisfies (\ref{O})
 and therefore, by Theorem \ref{un},
we have uniqueness for (\ref{R'}) with the given function 
$t(s)$, and we conclude that $x(s)=X(t(0),t(s),x(0)).$
\qed

\begin{cor}\label{cor}
Under the hypotheses of Theorem \ref{sup}, if in addition 
the vector-field $V$ satisfies (\ref{O}),
then $\mu_T=(X_0^T)_{\sharp}\mu_0$.
\end{cor}
\proof
Let $\mL\subset C([0,1],\Rm^{d+1})$ be a Borel set 
of full $\nu$-measure
formed by Lipschitz curves which satisfy all the properties 
of \ref{curve} in Theorem \ref{sup}.
If $\gamma(s)=(t(s),x(s))$ is a curve in $\mL$,
then $\gamma$ is one-to-one and, by Lemma \ref{su}, it
is of the form $\gamma(s)=(t(s),x(t(s))$.
We conclude that $t(s)$ is one to one and thus monotone.
By \ref{boundary} in Theorem \ref{sup} we can assume
in addition that
$t(0)\in \{0,T\}$ and $t(1)\in \{0,T\}$.

Denoting by $\mL^+$ the Borel subset of $\mL$
formed by curves $\gamma=(t,x)$ such that 
$t$ is increasing on $[0,1]$ and satisfies $t(0)=0$
and $t(1)=S$, and by $\mL^-$
the Borel subset of $\mL$
formed by curves $\gamma=(t,x)$ such that 
$t$ is decreasing on $[0,1]$ and satisfies $t(0)=T$
and $t(1)=0$,
we conclude that 
$\mL^+\cup\mL^-=\mL$.
We  denote by $\nu^\pm$ the restrictions of $\nu$
to $\mL^\pm$.
The measures $\nu^\pm$ are mutually singular, non-negative,
and $\nu=\nu^++\nu^-$.
Let 
$$B_i:\mL^+\cup \mL^-\lto \Rm^d
$$
be the Borel map defined by
$B_i(\gamma)=x(0)$ if $\gamma\in \mL^+$
and 
$B_i(\gamma)=x(1)$ if $\gamma\in \mL^-$.
Similarly, we define 
$$B_f:\mL^+\cup \mL^-\lto \Rm^d
$$
by
$B_i(\gamma)=x(0)$ if $\gamma\in \mL^-$
and 
$B_i(\gamma)=x(1)$ if $\gamma\in \mL^+$.
Note that 
$$B_f=X_0^T\circ B_i$$
on $\mL$.
We have the identities
$(ev_1)_{\sharp}\nu^+=\delta_T\otimes (B_f)_{\sharp}\nu^+,
(ev_0)_{\sharp}\nu^+=\delta_0\otimes (B_i)_{\sharp}\nu^+,
(ev_1)_{\sharp}\nu^-=\delta_0\otimes (B_i)_{\sharp}\nu^-$, and
$
(ev_0)_{\sharp}\nu^-=\delta_T\otimes (B_f)_{\sharp}\nu^-.
$
It follows that 
$$
\delta_T\otimes \mu_T
-\delta_0\otimes \mu_0=
(ev_1)_{\sharp}(\nu^++\nu^-)-(ev_0)_{\sharp}(\nu^++\nu^-)
=\delta_T\otimes (B_f)_{\sharp}(\nu^+-\nu^-)+
\delta_0\otimes (B_i)_{\sharp}(\nu^+-\nu^-).
$$
We conclude 
that $\mu_0=(B_i)_{\sharp}(\nu^--\nu^+)$
and  $\mu_S=(B_f)_{\sharp}(\nu^--\nu^+)$.
As a consequence, we have 
$$\mu_T=(X_0^T)_{\sharp}\mu_0.$$

\qed
\textsc{Proof of Theorem \ref{os} : }
We want to prove that $\mu_t=(X_s^t)_{\sharp}\mu_s$
for each $s$ and $t$ in $[0,T]$.
Since $X_t^s=(X_s^t)^{-1}$, it is enough to prove the statement
when $s<t$.
In order to do so, it is enough to apply Corollary \ref{cor}
on the time interval $[s,t]$ instead of $[0,T]$.
\qed
\appendix

\section{Topology on the spaces of measures}\label{top}
Let us first define the separable Banach space $C_0(\Rm^d)$
formed by the continuous functions which converge to zero at
infinity, with the uniform norm.
The space $\mM(\Rm^d)$ of bounded Borel measures
can be identified with the topological dual of $C_0(\Rm^d)$,
and the dual norm is $\|\mu\|=|\mu(\Rm^d)|$.
We endow this dual space with the weak-$*$ topology,
it is known that the unit ball $\mB:=\{\mu\in \mM | \|\mu\|\leq 1\}$
is compact and metrizable.
It is useful to work with  a specific distance.
In order to define this distance, we consider a sequence 
 $u_n$  of compactly supported smooth functions
on $\Rm^d$ which generates a dense vector subspace in $C_0(\Rm^d)$.
We assume in addition that $\|u_n\|_{C^1}=1$.
A distance on $\mB$ can be defined by the formula
$$
d(\mu,\eta)=\sum_{n\in \Nm} 
\frac{\big|\int u_n d\mu-\int u_n d\eta\big|}
{ 2^n}.
$$
It is well-known that the topology associated to this distance is
the weak-$*$ topology. $(\mB,d)$ is a compact metric space.
Now let $\mu_t$ satisfy (\ref{D}) for each of the functions $u_n$.
Then, given two times $s\leq t$, we have
$$\int u_n d\mu_t-\int u_nd\mu_s=\int_s^t\int_{\Rm^d} 
d(u_n)_x(V(\sigma,x))d\mu_{\sigma} d\sigma
$$
and therefore
\begin{equation}\tag{E}\label{E}
d(\mu_s,\mu_t)\leq \int_s^t \|V_{\sigma}\|_{\infty}d\sigma.
\end{equation}
We deduce that, given a Borel vector-field $V$ satisfying
(\ref{B}), the set of all solutions of (\ref{PDE}) in 
$C([0,T],\mB)$ is equi-continuous.

In order to put a topology on the space of measure-valued
solutions of (\ref{PDE}), a small digression is needed.
We define $L^1([0,T],C_0)$ as the set of Borel 
maps $u:[0,T]\lto C_0(\Rm^d)$ such that 
$\int_0^T \|u(t)\|_{\infty}dt<\infty.$
This is a separable banach space, which coincides with the space
of Borel measurable functions $v(t,x):[0,T]\times \Rm^d\lto \Rm$
such that, for each $t$, the map $v_t$ belongs to $C_0(\Rm^d)$
and such that, in addition, $\int_0^T \|v_t\|_{\infty}dt$ is finite.
All this is classical, see for example \cite{note}.
We denote by $\mY$ the unit ball in the dual of 
$L^1([0,T],C_0)$ endowed with the weak-$*$ topology.
This dual can be naturally identified with 
$L^{\infty}([0,T], \mM)$, and so it can be seen as a space of
Young measures, hence the notation.
$\mY$ is a compact metrizable space.
There is a natural embedding of $C([0,T],\mB)$
into $\mY$.
If $\mu^n_t$ is an equi-continuous sequence 
in $C([0,T],\mB)$ which converges in $\mY$ 
to the limit $\mu_t\in \mY$, then 
$\mu_t\in C([0,T],\mB)$ and the convergence is uniform in 
$C([0,T],\mB)$.

\section{Set valued maps}\label{set}
The classical reference for this material is the book of Kuratowski, \cite{K}. Let X be 
a complete metric space, and K be a compact metric space. A set-valued map S 
associates to each point $x \in X$ a subset $S (x)$ of $K$. 
The set-valued map S is called 
upper semi-continuous if its graph 
$$\{(x, y) \in  X \times  K \quad| \quad y \in S (x)\} 
$$
is closed. We consider from now on an upper semi-continuous set function S. 
Given $U \subset  K$, we define $S^{-1} (U )$
 as the set of points $x \in B$ such that $S (x) \subset  U . $
It is easy to see that $S^{-1} (U )$  is open 
for each open set $U$ ( recall that $S$ is 
upper semi-continuous). Since every closed 
set is a $G_{\delta}$ (a countable intersection of open sets), we get: 
\begin{lem} 
If $S$ is upper semi-continuous, then $S^{ -1} (F )$ 
is a $G_{\delta}$ for each closed 
set $F \subset K$.
\end{lem} 
We say that $x$ is a point of continuity of $S$ if, for each 
$y \in S (x)$ and each sequence 
$x_n \lto x$ in $X$, there exists a sequence 
$y_n \lto y$ such that $y_n \in S (x_n ).$ 

\begin{thm}\label{K}
 If S is an upper semi-continuous set function, then the set of points 
of continuity of S is a dense $G_{\delta}$. 
\end{thm}
\proof
 Let $U_k$ be a countable base of open sets, and let $F_k$ be the complement of 
$U_k$ . 
We claim that the set of points of continuity is 
$$
\bigcap_{k\in N}
\left[
\big(\overline{S^{-1} (F_k ) }\big)^c
\cup S^{-1} (F_k )
 \right].
 $$
Each of the sets 
$$
\big(\overline{S^{1} (F_k ) }\big)^c
\cup S^{1} (F_k )$$
is a  $G_{\delta}$ because $S^{-1} (F_k )$ is a  $G_{\delta}$
and 
$
\big(\overline{S^{-1} (F_k ) }\big)^c
$
is open.
In addition, it is clearly 
dense. By the Baire property, we conclude that the set of continuity points is a 
dense $G_{\delta}$.
 We now have to check the claim. The point $x$ is not a point of continuity 
if and only if there exists an open set $U_k$
 such that $ U_k\cap S (x)$ is not empty and 
a sequence $x_n \lto x$  such that $ S (x_n ) \cap U_k$ is empty.
 This amounts to say that 
 $x\not \in  S^{-1} (F_k )$ an $x\in \overline{S^{-1} (F_k ) }$.
As a consequence, the complement of the set of 
continuity points is 
$$
\bigcup_k \big( \overline{S^{-1} (F_k ) }-  S^{-1} (F_k ).\big)
$$
\qed

\end{document}